\documentclass{amsproc}
\usepackage{amssymb}
\usepackage{amsfonts}

\theoremstyle{plain}

\newtheorem{definition}{Definition}

\newtheorem{lemma}{Lemma}

\newtheorem{theorem}{Theorem}
\numberwithin{equation}{section}
\input{tcilatex}

\begin{document}
\title[TOEPLITZNESS OF COMPOSITION OPERATORS]{DECAY ESTIMATES CONNECTED WITH
TOEPLITZNESS OF COMPOSITION OPERATORS}
\author{Faruk F. Abi-Khuzam}
\address{Department of Mathematics, American University of Beirut, Beirut
Lebanon}
\email{farukakh@aub.edu.lb}
\date{July 10, 2011}

\begin{abstract}
Let $f\in L^{\infty }(T^{d})$ with $\Vert f\Vert _{L^{\infty }(T^{d})}\leq 1$%
, $\nu \in 
\mathbb{Z}
^{d},n,k\in 
\mathbb{Z}
$ and put $b_{n,n-k}=\int_{E}f(x)^{n}e^{-2\pi i(n-k)\nu \cdot x}dx$, $%
E=\{x\in T^{d}:|f(x)|=1\}.$ Shayya conjectured that, if $\hat{f}(\xi )=0$
for all $\xi $ in a half-space $S$ of lattice points, and $\nu \in -S$, and $%
\hat{f}(0)\neq 0$, then $\lim_{n\rightarrow \infty }b_{n,n-k}=0,k\in 
\mathbb{Z}
$. This is a higher dimensional version of an earlier conjecture of Nazarov
and Shapiro, the truth of which would imply that any composition operator is
weakly asymptotically Toeplitz on the Hardy space $H^{2}$. Shayya proved
that the arithmetic means of the main diagonal entries, $\{b_{n,n}\}$, decay
like $\{\log N\}^{-1}$, but the proof does not extend to other diagonal
entries. Using an elementary method, we extend and improve all known
results. In particular, we show that the arithmetic means of the entries of
each diagonal, $\{b_{n,n-k}\},$ decay like $\{\log N\log \log N\}^{-1}$
uniformly in $k\in 
\mathbb{Z}
$.
\end{abstract}

\maketitle

\section{INTRODUCTION}

Let $\phi $ be a function holomorphic in the open unit disk $U$ of the
complex plane and satisfying the two conditions%
\begin{equation*}
\phi (0)\neq 0,|\phi (z)|<1,z\in U.
\end{equation*}%
The composition operator $C_{\phi }:H(U)\rightarrow H(U)$, defined via $\phi 
$ by $C_{\phi }f=f\circ \phi $, restricts to a bounded linear operator on
the Hardy space $H^{2}(U)$, and so has an infinite matrix $(c_{\alpha .\beta
})$ associated to it and defined by 
\begin{equation*}
c_{\alpha ,\beta }=\int_{\partial U}\phi (z)^{\alpha }z^{-\beta }\frac{dz}{%
2\pi iz}.
\end{equation*}%
If $S$ is the forward shift operator on $H^{2}(U)$, and $S^{\ast }$ is its
adjoint, then $C_{\phi }$ is said to be weakly asymptotically Toeplitz , or
WAT, if the sequence $\{S^{\ast n}C_{\phi }S^{n}\}$ of bounded linear
operators on $H^{2}(U)$ converges in the weak operator topology.
Equivalently, $C_{\phi }$ is WAT, if 
\begin{equation*}
\lim_{n\rightarrow \infty }c_{n,n-k}=0,k\in 
\mathbb{Z}
\text{.}
\end{equation*}%
If $E=\{z\in \partial U:|\phi (z)|=1\}$, and we put 
\begin{equation*}
b_{n,n-k}=\int_{E}\phi (z)^{n}z^{k-n}\frac{dz}{2\pi iz},
\end{equation*}%
then, clearly, $c_{n,n-k}\rightarrow 0$, if and only if $b_{n,n-k}%
\rightarrow 0$, as $n\rightarrow \infty .$The WAT conjecture on the circle,
formulated by Nazarov and Shapiro [4], is the statement that 
\begin{equation}
\lim_{n\rightarrow \infty }b_{n,n-k}=0,k\in 
\mathbb{Z}
.
\end{equation}%
At present, it is not known if \ $(1.1)$ is true, even when $k=0$. But the
arithmetic means are more tame. In [5], it was shown that 
\begin{equation}
\lim_{N\rightarrow \infty }\frac{1}{N}\sum_{n=1}^{N}|b_{n,n-k}|^{2}=0,k\in 
\mathbb{Z}
\text{.}
\end{equation}%
An infinite series carrying the full sequence $\{c_{n,n-k}\}$ was considered
in [1], where it was shown, using complex analytic methods, that 
\begin{equation*}
\sum_{n=1}^{\infty }\frac{|c_{n,n-k}|^{2}}{n}<\infty ,k\geq 0,
\end{equation*}%
from which it follows, through a simple summation by parts argument, that
the arithmetic means of the sequence $\{|c_{n,n-k}|^{2}\}$ tend to $0$, and
also that 
\begin{equation*}
\sum_{n=1}^{\infty }\frac{1}{n(n+1)}\sum_{j=1}^{n}|c_{j,j-k}|^{2}<\infty
,k\geq 0.
\end{equation*}%
The convergence of this later series implies that \ $\liminf_{n\rightarrow
\infty }\frac{\log n}{n}\sum_{j=1}^{n}|c_{j,j-k}|^{2}=0$, and this
anticipates later work on the arithmetic means of the sequence $%
\{|b_{n,n-k}|^{2}\}$. The convergence of $\sum_{n=1}^{\infty }\frac{%
|b_{n,n-k}|^{2}}{n}$ was established as part of a study of these problems
undertaken recently by Shayya [6] in the setting of functions defined on the 
$d$-dimensional torus $T^{d}$. He proved the existence of a constant $C$,
essentially dependent on $\phi (0)$, such that 
\begin{equation}
\frac{1}{N}\sum_{n=M}^{M+N}|b_{n,n}|^{2}\leq \frac{C}{\log \frac{N}{\pi }}%
,M\geq 1,N\geq 2,
\end{equation}%
thereby obtaining a quantitative strengthening of $(1.2)$ , in the case $k=0$%
. Shayya's method uses a positive Borel measure $\mu $ such that 
\begin{equation*}
|b_{n,n}|^{2}=\int_{-1/2}^{1/2}e^{-2\pi in\theta }d\mu (\theta ),
\end{equation*}%
along with the convergence of the series $\sum_{n=1}^{\infty }\frac{%
|b_{n,n}|^{2}}{n}$. This measure, obtained as an application of Herglotz's
theorem, is only available when $k=0$.

Shayya's result $(1.3)$ tells us that the $\lim \sup_{n\rightarrow \infty }%
\frac{\log n}{n}\sum_{j=1}^{n}|b_{n,n}|^{2}<\infty $, while, by the argument
above, $\lim \inf_{n\rightarrow \infty }\frac{\log n}{n}%
\sum_{j=1}^{n}|b_{n,n}|^{2}=0$. So if the limit exists, it would have to be $%
0$, and a natural question would then be to find the rate of decay to $0$.
In the present work we consider the more general problem of the arithmetic
means of $|b_{n,n-k}|^{2},k\in 
\mathbb{Z}
.$ We obtain, uniformly in $k\in 
\mathbb{Z}
$, a quantitative strengthening of $(1.2)$. In addition, we show that these
means actually decay much faster than suggested by $(1.3)$ and we obtain,
among other results, the inequality 
\begin{equation}
\frac{1}{N+3}\sum_{n=M}^{M+N}|b_{n,n-k}|^{2}\leq \frac{C/(1-1/\log 3)}{\log
(N+3)\log [\frac{\log (N+4)}{\log 4}]},k\in 
\mathbb{Z}
,
\end{equation}%
for all $M\geq 1$, and $N\geq 1$, where the constant $C$ depends only on $%
\phi (0)$.

\bigskip

\section{PRELIMINARIES AND STATEMENT OF RESULTS}

We shall consider the problem of the decay of the arithmetic means in the
setting of functions defined on $T^{d}$, as formulated in [6]. Since we
shall be working on the coefficients $|b_{n,n-k}|^{2}$ with $k\in 
\mathbb{Z}
,$ we will not have available a measure as in [6], and so we shall proceed
along different, rather elementary, lines. However, there will be two
important common ingredients between the present paper and [6], namely the
use of Helson's and Lowdenslager's generalization of Szego's theorem [3],
and the restriction to the set where the boundary values of the given
function have unit modulus. We shall obtain a bound on the sum of the series 
\begin{equation*}
\sum_{\substack{ m=-\infty  \\ m\neq N}}^{\infty }\frac{|b_{m,m-k}|^{2}}{%
|m-N|},
\end{equation*}%
where $k,N\in 
\mathbb{Z}
$. The usefulness of this stems from the fact that the bound will be
independent of $k$ and $N$, thereby making it possible to obtain uniform
bounds on the arithmetic means of \ $\{|b_{m,m-k}|^{2}\}$.

We now give a summary, taken mostly from [6], of material needed in
connection with the Helson-Lowdenslager result.

\begin{definition}
A set $S\subset 
\mathbb{Z}
^{d}$ is said to be a half-space of lattice points in $%
\mathbb{R}
^{d}$ if
\end{definition}

\textit{(i) }$0\not\in S$

\textit{(ii) }$\xi \neq 0$\textit{\ implies }$\xi \in S$\textit{\ or }$-\xi
\in S$

\textit{(iii) }$\xi ,\xi ^{\prime }\in S$\textit{\ implies }$\xi +\xi
^{\prime }\in S.$

Clearly, if $S$ is a half-space of lattice points in $%
\mathbb{R}
^{d}$, then so is its reflection in the origin. i.e. $-S=\{-\xi :\xi \in
S\}. $

Also [6], if $S_{1}$ and $S_{2}$ are half-spaces of lattice points in $%
\mathbb{R}
^{d_{1}}$ and $%
\mathbb{R}
^{d_{2}}$ respectively, then $S_{1}\times 
\mathbb{R}
^{d_{2}}\cup 
\mathbb{R}
^{d_{1}}\times S_{2}$ contains a half-space of lattice points in $%
\mathbb{R}
^{d_{1}}\times 
\mathbb{R}
^{d_{2}}$(e.g., $S_{1}\times 
\mathbb{R}
^{d_{2}}\cup \{0\}\times S_{2}$).

\bigskip We shall need the following generalization of Szego's theorem [3],
stated as corollary A in [6].

\begin{lemma}
Suppose $f\in L^{2}(T^{d})$, $S$ is a half-space of lattice points in $%
\mathbb{R}
^{d}$, and $\hat{f}(\xi )=0$ for all $\xi \in S$. then 
\begin{equation*}
\int_{T^{d}}\log |f(x)|dx\geq \log |\hat{f}(0)|.
\end{equation*}
\end{lemma}

As a simple application of this lemma, suppose $S$ is a half-space of
lattice points in $%
\mathbb{R}
^{d}$, $f\in L^{\infty }(T^{d}),\Vert f\Vert _{L^{\infty }(T^{d})}\leq 1$,
and $\hat{f}(\xi )=0$ for all $\xi \in S$. For $0<r<1$, and $\nu \in -S$,
define $F:T^{d}\times T^{d}\rightarrow 
\mathbb{C}
$ by 
\begin{equation*}
F(x,y)=e^{2\pi i(\nu ,-\nu )\cdot (x,y)}-rf(x)\overline{f(y)}.
\end{equation*}%
Then, [6], $\hat{F}(\xi ,\eta )=0$ for all $(\xi ,\eta )\in S\times 
\mathbb{R}
^{d}\cup 
\mathbb{R}
^{d}\times (-S)$, and so for all $(\xi ,\eta )$ in a half-space of lattice
points in $%
\mathbb{R}
^{d}\times 
\mathbb{R}
^{d}$. Lemma 1 applies and we conclude that 
\begin{equation*}
\int_{T^{d}}\int_{T^{d}}\log |F(x,y)|dxdy\geq \log |\hat{F}(0,0)|=\log r|%
\hat{f}(0)|^{2}.
\end{equation*}

Also, using the equality $|\log |u||=2\log ^{+}|u|-\log |u|$, we have, for
any measurable subset of $T^{d}$, 
\begin{equation}
\int_{E}\int_{E}|\log |F(x,y)||dxdy\leq \int_{T^{d}}\int_{T^{d}}|\log
|F(x,y)||dxdy
\end{equation}%
\begin{equation*}
=2\int_{T^{d}}\int_{T^{d}}\log ^{+}|F(x,y)|dxdy-\int_{T^{d}}\int_{T^{d}}\log
|F(x,y)|dxdy
\end{equation*}%
\begin{equation*}
\leq 2\log (1+r)-\log r|\hat{f}(0)|^{2}\leq \log \frac{4}{r|\hat{f}(0)|^{2}}.
\end{equation*}

We shall also need two sequences of functions $\log _{j}x$ and $L_{q}(x)$.
The first is the iterated logarithmic function, and the second is a product
of different iterates:

If $j\geq 1$ is an integer, we define the function $\log _{j}$ as follows:%
\begin{equation*}
\log _{1}x=\log x,\log _{j}x=\log \circ \log _{j-1}x\text{.}
\end{equation*}

We will need to use the iterated logarithm on the domain where it is
well-defined and positive, and this, of course, is trivial to determine.In
terms of the iterated logarithm we define the function $L_{q}$ by 
\begin{equation*}
L_{q}(x)=\dprod\limits_{j=1}^{q}\log _{j}x,q\geq 1.
\end{equation*}%
Also, for a positive differentiable function $g$, we define [2] the function 
$a(.;g)$ by $a(x;g)=xg^{\prime }(x)/g(x)$, and note that 
\begin{eqnarray*}
a(x;L_{1}) &=&\frac{1}{\log x},a(x;L_{q})=\frac{1}{\log x}\left(
1+\sum_{j=2}^{q}\frac{1}{\log _{j}x}\right) ,q\geq 2. \\
\log _{q+1}^{\prime }(x) &=&\frac{1}{xL_{q}(x)}.
\end{eqnarray*}%
Clearly, $a(x;L_{q})$ decreases monotonically to $0$ as $x\rightarrow \infty 
$, and so there exist positive constants $\alpha _{q}$ and $\gamma _{q}$ ,
depending only on the integer $q$, such that $0<\alpha _{q}<1$, and 
\begin{equation*}
\log _{q+1}x>0,0<a(x;L_{q})<\alpha _{q},x\geq \gamma _{q}\text{.}
\end{equation*}

\begin{theorem}
Suppose $f\in L^{\infty }(T^{d})$ with $\Vert f\Vert _{L^{\infty
}(T^{d})}\leq 1,$ and $\nu \in 
\mathbb{Z}
^{d}$. For $n,k\in 
\mathbb{Z}
$, define 
\begin{equation*}
b_{n,n-k}=\int_{E}f(x)^{n}e^{-2\pi i(n-k)\nu \cdot x}dx,
\end{equation*}%
where $E=\{x\in T^{d}:|f(x)|=1\}.$When $\hat{f}(0)\neq 0$, \textit{put }$%
C=\log \frac{16}{|\hat{f}(0)|^{4}}$\textit{.}
\end{theorem}

\textit{(i) If\ }$\hat{f}(\xi )=0$\textit{\ for all }$\xi $\textit{\ in a
half-space }$S$\textit{\ of lattice points, }$\nu \in -S$\textit{, }$\hat{f}%
(0)\neq 0$, \textit{and }$N\in 
\mathbb{Z}
$\textit{, then }%
\begin{equation}
\sum_{\substack{ m=-\infty  \\ m\neq N}}^{\infty }\frac{|b_{m,m-k}|^{2}}{%
|m-N|}\leq \log \frac{16}{|\hat{f}(0)|^{4}},k\in 
\mathbb{Z}
.
\end{equation}

\textit{(ii) If\ }$\hat{f}(\xi )=0$\textit{\ for all }$\xi $\textit{\ in a
half-space }$S$\textit{\ of lattice points, }$\nu \in -S$\textit{, }$\hat{f}%
(0)\neq 0$\textit{, }$M\geq 1$\textit{, and }$p\geq 1$\textit{, then}%
\begin{equation}
\frac{1}{p+1}\sum_{m=M}^{M+p}|b_{m,m-k}|^{2}\leq \frac{C}{\log (p+1)},k\in 
\mathbb{Z}
.
\end{equation}

\textit{(iii) Suppose }$g$ \textit{is positive, and }$0<a(x;g)<\alpha $%
\textit{, for all }$x\geq \gamma \geq 1$. \textit{If\ }$\hat{f}(\xi )=0$%
\textit{\ for all }$\xi $\textit{\ in a half-space }$S$\textit{\ of lattice
points, }$\nu \in -S$\textit{, }$\hat{f}(0)\neq 0$\textit{, }$M\geq 1$%
\textit{, and }$p\geq 1$\textit{, then }%
\begin{equation}
\left( \int_{1}^{p+1}\frac{dt}{tg(t+\gamma )}\right)
\sum_{m=M}^{M+p}|b_{m,m-k}|^{2}\leq \frac{C}{1-\alpha }\cdot \frac{p+\gamma 
}{g(p+\gamma )},k\in 
\mathbb{Z}
.
\end{equation}

\textit{(iv) Suppose }$L_{q},\alpha _{q}$\textit{, and }$\gamma _{q}$\textit{%
\ are as defined above.} \textit{If\ }$\hat{f}(\xi )=0$\textit{\ for all }$%
\xi $\textit{\ in a half-space }$S$\textit{\ of lattice points, }$\nu \in -S$%
\textit{, }$\hat{f}(0)\neq 0$\textit{, }$M\geq 1$\textit{, and }$p\geq 1$%
\textit{, then}%
\begin{equation}
\frac{1}{p+\gamma _{q}}\sum_{m=M}^{M+p}|b_{m,m-k}|^{2}<\frac{C/(1-\alpha
_{q})}{L_{q}(p+\gamma _{q})[\log _{q+1}(p+1+\gamma _{q})-\log
_{q+1}(1+\gamma _{q})]},k\in 
\mathbb{Z}
\text{.}
\end{equation}%
\textit{Also, for a fixed }$q\geq 1$\textit{, there is a constant }$%
C^{\prime }$\textit{\ depending only on }$q$ and $C$\textit{, such that, for
all large }$p$\textit{,}%
\begin{equation*}
\frac{1}{p}\sum_{m=M}^{M+p}|b_{m,m-k}|^{2}<\frac{C^{\prime }}{L_{q+1}(p)}=%
\frac{C^{\prime }}{\log p\cdot \log \log p\cdot \cdot \cdot \log _{q+1}p},
\end{equation*}%
\textit{uniformly in }$k\in 
\mathbb{Z}
$\textit{.}

The inequality in part (ii) represents a sharpening, and an extension of $%
(1.3)$. The inequality in part (iv) shows that, uniformly in $k\in 
\mathbb{Z}
$, the decay of the arithmetic means of the sequence $\{|b_{m,m-k}|^{2}\}$
is faster than that of $L_{q+1}^{-1}$ , for any $q\geq 1.$

\section{PROOF OF THEOREM 1}

The definition of $b_{n,n-k}$ gives us that \textit{\ }%
\begin{equation*}
|b_{n.n-k}|^{2}=\int_{E}\int_{E}\left( f(x)\overline{f(y)}e^{-2\pi i\nu
\cdot (x-y)}\right) ^{n}.e^{2\pi ik\nu \cdot (x-y)}dxdy,
\end{equation*}%
and, for $N,k\in 
\mathbb{Z}
$, we have, since $|f(x)|=1$ for $x\in E,$\textit{\ }%
\begin{equation*}
|b_{n+N,n+N-k}|^{2}+|b_{-n+N,-n+N-k}|^{2}
\end{equation*}%
\begin{equation*}
=2\int_{E}\int_{E}\left( f(x)\overline{f(y)}e^{-2\pi i\nu \cdot
(x-y)}\right) ^{N}.e^{2\pi ik\nu \cdot (x-y)}\func{Re}\left( f(x)\overline{%
f(y)}e^{-2\pi i\nu \cdot (x-y)}\right) ^{n}dxdy.
\end{equation*}%
Multiplying by $r^{n},0<r<1$, dividing by $n$ and summing over $n\geq 1,$ we
obtain \textit{\ }%
\begin{equation*}
\sum_{n=1}^{\infty }\frac{|b_{n+N,n+N-k}|^{2}+|b_{-n+N,-n+N-k}|^{2}}{n}r^{n}
\end{equation*}%
\begin{equation*}
=2\int_{E}\int_{E}\left( f(x)\overline{f(y)}e^{-2\pi i\nu \cdot
(x-y)}\right) ^{N}.e^{2\pi ik\nu \cdot (x-y)}.\log \frac{1}{|F(x,y)|}dxdy,
\end{equation*}%
where $\ F(x,y)=e^{2\pi i\nu \cdot (x-y)}-rf(x)\overline{f(y)}\ $. In view
of the estimate $(2.1)$ on the integral of $|\log |F(x,y)||$ we get that 
\textit{\ }%
\begin{equation*}
\sum_{n=1}^{\infty }\frac{|b_{n+N,n+N-k}|^{2}+|b_{-n+N,-n+N-k}|^{2}}{n}%
r^{n}\leq \log \frac{16}{r^{2}|\hat{f}(0)|^{4}},
\end{equation*}%
and letting $r\rightarrow 1-$, we obtain\textit{\ }%
\begin{equation*}
\sum_{n=1}^{\infty }\frac{|b_{n+N,n+N-k}|^{2}+|b_{-n+N,-n+N-k}|^{2}}{n}\leq
C.
\end{equation*}%
Reindexing we arrive at\textit{\ }%
\begin{equation*}
\sum_{\substack{ m\neq N  \\ m=-\infty }}^{\infty }\frac{|b_{m,m-k}|^{2}}{%
|m-N|}\leq C,k\in 
\mathbb{Z}
,
\end{equation*}%
and this completes the proof of part (i).

If $M\geq 1,p\geq 1$, then, for a given $N$, the partial sums in $(2.2)$ for 
$M\leq m\leq M+p$, $m\neq N$, are also bounded by the constant $C$ uniformly
in $k,N\in 
\mathbb{Z}
$, and we may sum them over $N$ from $M$ to $M+p$ to get, after an
interchange in the order of summation

\begin{equation}
\sum_{m=M}^{M+p}\sum_{\substack{ N=M  \\ N\neq m}}^{M+p}\frac{|b_{m,m-k}|^{2}%
}{|m-N|}\leq \sum_{N=M}^{M+p}C=C(p+1),k\in 
\mathbb{Z}
.
\end{equation}%
On the left-hand side of $(3.1)$ we use the inequality \textit{\ }$\sum 
_{\substack{ N=M  \\ N\neq m}}^{M+p}\frac{1}{|m-N|}\geq \sum_{j=1}^{p}\frac{1%
}{j}$ :\textit{\ }%
\begin{equation*}
\left( \sum_{m=M}^{M+p}|b_{m},_{m-k}|^{2}\right) \left( \sum_{\substack{ N=M 
\\ N\neq m}}^{M+p}\frac{1}{|m-N|}\right)
\end{equation*}%
\begin{equation*}
\geq \left( \sum_{m=M}^{M+p}|b_{m},_{m-k}|^{2}\right) \left( \sum_{j=1}^{p}%
\frac{1}{j}\right) \geq \log (p+1)\left(
\sum_{m=M}^{M+p}|b_{m},_{m-k}|^{2}\right) ,
\end{equation*}%
so that \textit{\ }%
\begin{equation*}
\log (p+1)\sum_{m=M}^{M+p}|b_{m},_{m-k}|^{2}\leq C(p+1),k\in 
\mathbb{Z}
,
\end{equation*}%
and the inequality in $(2.3)$ follows.

For proof of part (iii), first divide both sides of $(2.2)$ by $g(N-M+\gamma
),N\geq M\geq 1,$and\ proceed exactly as in the proof of part (ii). In the
resulting inequality the left-hand side is\textit{\ }%
\begin{equation*}
(\sum_{m=M}^{M+p}|b_{m},_{m-k}|^{2})\sum_{\substack{ N=M  \\ N\neq m}}^{M+p}%
\frac{1}{|m-N|g(N-M+\gamma )}\geq
(\sum_{m=M}^{M+p}|b_{m},_{m-k}|^{2})\sum_{j=1}^{p}\frac{1}{jg(j+\gamma )}
\end{equation*}%
\begin{equation*}
\geq \sum_{m=M}^{M+p}|b_{m},_{m-k}|^{2}\int_{1}^{p+1}\frac{dt}{tg(t+\gamma )}%
,
\end{equation*}%
while the right-hand side is\textit{\ }%
\begin{equation*}
C\sum_{N=M}^{M+p}\frac{1}{g(N-M+\gamma )}<C\{\frac{1}{g(\gamma )}%
+\int_{0}^{p}\frac{dt}{g(t+\gamma )}\}.
\end{equation*}%
The conditions on\textit{\ }$a(x;g),$ along with the positivity of $g$,
imply that both functions $g$ and $x/g(x)$ are increasing on $[\gamma
,\infty )$. An application of the Cauchy mean-value theorem gives a $\xi \in
(\gamma ,x)$ such that\textit{\ }%
\begin{equation*}
\frac{\int_{\gamma }^{x}\frac{dt}{g(t)}}{\frac{x}{g(x)}-\frac{\gamma }{%
g(\gamma )}}=\frac{1}{1-a(\xi ;g)}<\frac{1}{1-\alpha }.
\end{equation*}%
Since\textit{\ }$\gamma \geq 1,0<\alpha <1,$we have%
\begin{equation*}
-\left( \frac{1}{1-\alpha }\right) \frac{\gamma }{g(\gamma )}+\frac{1}{%
g(\gamma )}<0,
\end{equation*}%
and it follows that\textit{\ }%
\begin{equation}
\frac{1}{g(\gamma )}+\int_{0}^{p}\frac{dt}{g(t+\gamma )}<\frac{1}{1-\alpha }%
\cdot \frac{p+\gamma }{g(p+\gamma )},
\end{equation}%
which gives the required bound and the inequality in $(2.4)$ follows.

\textit{If }$q\geq 1,$ and\textit{\ }$L_{q}$ is the function defined in
section 2, and $\alpha _{q}$ and $\gamma _{q}$ are the positive constants
associated with it, we may apply \ the result in part (iii) with $%
g(x)=L_{q}(x)$. This gives us%
\begin{equation*}
\left( \int_{1}^{p+1}\frac{dt}{tL_{q}(t+\gamma _{q})}\right)
\sum_{m=M}^{M+p}|b_{m,m-k}|^{2}\leq \frac{C}{1-\alpha _{q}}\cdot \frac{%
p+\gamma _{q}}{L_{q}(p+\gamma _{q})},k\in 
\mathbb{Z}
.
\end{equation*}%
For the integral on the left-hand side we have \textit{\ }%
\begin{equation*}
\int_{1}^{p+1}\frac{dt}{tL_{q}(t+\gamma _{q})}\geq \int_{1}^{p+1}\frac{dt}{%
(t+\gamma _{q})L_{q}(t+\gamma _{q})}
\end{equation*}%
\begin{equation*}
=\log _{q+1}(p+1+\gamma _{q})-\log _{q+1}(1+\gamma _{q}),
\end{equation*}%
and it follows that\textit{\ }%
\begin{equation*}
\frac{1}{p+\gamma _{q}}\sum_{m=M}^{M+p}|b_{m,m-k}|^{2}\leq \frac{C/(1-\alpha
_{q})}{L_{q}(p+\gamma _{q})[\log _{q+1}(p+1+\gamma _{q})-\log
_{q+1}(1+\gamma _{q})]},
\end{equation*}%
for all $k\in 
\mathbb{Z}
$.

If $q\geq 1$ is fixed, since $L_{q}(p+\gamma _{q})[\log _{q+1}(p+1+\gamma
_{q})-\log _{q+1}(1+\gamma _{q})]\thicksim L_{q+1}(p)$ as \textit{\ }$%
p\rightarrow \infty ,$ there is a constant $C^{\prime }$depending only on $q$
and $C,$ such that, for all large $p$, \textit{\ }%
\begin{equation*}
\frac{1}{p}\sum_{m=M}^{M+p}|b_{m,m-k}|^{2}<\frac{C^{\prime }}{L_{q+1}(p)},
\end{equation*}%
uniformly in $k\in 
\mathbb{Z}
$.

In the special case where $q=1$, $L_{q}(x)=\log x$ and we may take $\gamma
=3 $, and $\alpha =\frac{1}{\log 3}.$ Then Theorem 1 gives us that \textit{\ 
}%
\begin{equation*}
\frac{1}{p+3}\sum_{m=M}^{M+p}|b_{m,m-k}|^{2}\leq \frac{C/(1-1/\log 3)}{\log
(p+3)[\log \log (p+4)-\log \log 4]},k\in 
\mathbb{Z}
\text{.}
\end{equation*}

There are some interesting questions that present themselves in connection
with the study of the problems considered in this note. One question
concerns the rate of decay of the tail end of the series $\sum_{n=1}^{\infty
}\frac{|b_{n,n-k}|^{2}}{n}$. Even in the case $k=0,$ this is not known.
Another question concerns the series $\sum_{n=1}^{\infty }L_{q}(n)\frac{%
|b_{n,n-k}|^{2}}{n}$. It is not known if this series is convergent, even in
the case $q=1$, and $k=0$.

\section{REFERENCES}

[1] F. ABI\_KHUZAM AND B. SHAYYA, \textit{A remark on the WAT conjecture,
preprint}, 2009.

[2]W. HAYMAN, \textit{A generalisation of Stirling's formula}, J. Reine.
Angew. Math. \textbf{196} (1956), 67-95.

[3]H. HELSON AND D. LOWDENSLAGER,\textit{\ Prediction theory and Fourier
series in several variables}, Acta Math. \textbf{99} (1958), 165-202.

[4] F. NAZAROV AND J. SHAPIRO, \textit{On the Toeplitzness of composition
operators}, Complex Var. Elliptic Wqu. \textbf{52} (2007), 193-210.

[5] J. SHAPIRO, \textit{Every composition operator is ( mean) asymptotically
Toeplitz}, J. Math. Anal. Appl. \textbf{333} (2007), 523-529.

[6] B.SHAYYA, \textit{The WAT conjecture on the Torus}, Proc. Amer. Math.
Soc. \textbf{139} ( 2011), 3633-3643.

\end{document}